\documentclass[a4paper,10pt]{article}
\usepackage{stmaryrd}
\usepackage{amsfonts}
\usepackage{bbm}
\usepackage{chemarrow}
\usepackage{amscd}
\usepackage{mathrsfs}
\usepackage{latexsym,amssymb,amsmath,amscd,amscd,amsthm,amsxtra}
\usepackage[dvips]{graphicx}
\usepackage[utf8]{inputenc}
\usepackage[T1]{fontenc}
\usepackage{lmodern}
\usepackage{amssymb}
\usepackage[all]{xy}
\usepackage{nicefrac,mathtools,enumitem}
\usepackage{microtype}
\usepackage{xcolor}
\newtheorem{thm}{Theorem}[section]
\newtheorem{lem}[thm]{Lemma}
\newtheorem{cor}[thm]{Corollary}
\newtheorem{pro}[thm]{Proposition}
\newtheorem{ex}[thm]{Example}
\newtheorem{rmk}[thm]{Remark}
\newtheorem{defi}[thm]{Definition}

\newcommand {\emptycomment}[1]{}

\newcommand {\nrn}[1]{   [    #1   ]^{\rm{3Lie}}}

\newcommand{\lon }{\,\rightarrow\,}
\newcommand{\be }{\begin{equation}}
\newcommand{\ee }{\end{equation}}

\newcommand{\pf}{\noindent{\bf Proof.}\ }

\newcommand{\g}{\mathfrak g}
\newcommand{\h}{\mathfrak h}

\newcommand{\huaF}{\mathcal{F}}

\newcommand{\huaU}{\mathcal{U}}

\newcommand{\huaX}{\mathcal{X}}
\newcommand{\huaY}{\mathcal{Y}}

\newcommand{\huaC}{{\mathcal{C}}}
\newcommand{\huaD}{\mathcal{D}}

\newcommand{\huaMC}{\mathcal{MC}}

\newcommand{\frke}{\mathfrak e}

\newcommand{\frkg}{\mathfrak g}
\newcommand{\frkh}{\mathfrak h}

\newcommand{\frkk}{\mathfrak k}

\newcommand{\frks}{\mathfrak s}

\newcommand{\frkX}{\mathfrak X}
\newcommand{\frkY}{\mathfrak Y}

\def\qed{\hfill ~\vrule height6pt width6pt depth0pt}

\newcommand{\Id}{\rm{Id}}

\newcommand{\br}[1]{   [ \cdot,    \cdot  ]   }

\newcommand{\dM}{\mathrm{d}}

\newcommand{\Hom}{\mathrm{Hom}}
\newcommand{\Der}{\mathrm{Der}}
\newcommand{\Lie}{\mathrm{Lie}}
\newcommand{\Rep}{\mathrm{Rep}}
\newcommand{\ab}{\mathrm{ab}}
\newcommand{\ob}{\mathrm{ob}}
\newcommand{\Set}{\mathrm{Set}}

\newcommand{\gl}{\mathfrak {gl}}

\newcommand{\ad}{\mathrm{ad}}
\newcommand{\pr}{\mathrm{pr}}

\newcommand{\K}{\mathbb{K}}

\newcommand{\T}{\mathbb{T}}

\begin{document}
\title{
{On non-abelian extensions of 3-Lie algebras
\thanks
 {
This research is supported by NSFC (11471139) and NSF of Jilin Province (20170101050JC).
 }
} }
\author{Lina Song$^1$, Abdenacer Makhlouf$^2$ and  Rong Tang$^1$
\\
\\
$^1$Department of Mathematics, Jilin University,
 Changchun 130012,  China \\
$^2$University of Haute Alsace, Laboratoire de Math\'ematiques,\\ \vspace{2mm}Informatique et Applications, Mulhouse, France\\\vspace{2mm}
 Email:
songln@jlu.edu.cn, abdenacer.makhlouf@uha.fr, tangrong16@mails.jlu.edu.cn
 }
\date{}
\maketitle

\begin{abstract}
In this paper, we study non-abelian extensions of 3-Lie algebras through Maurer-Cartan elements. We show that there is a one-to-one correspondence between isomorphism classes of non-abelian extensions of 3-Lie algebras and equivalence classes of Maurer-Cartan elements in a DGLA. The structure of the Leibniz algebra on the space of fundamental objects is also analyzed.
\end{abstract}


\section{Introduction}
Ternary Lie algebras (3-Lie algebras) or more generally $n$-ary Lie algebras are a natural generalization of Lie  algebras.  They were introduced and studied first by  Filippov in  \cite{Filippov}. This type of algebras appeared also in the algebraic formulation of Nambu Mechanics \cite{N}, generalizing Hamiltonian mechanics by considering two hamiltonians, see \cite{T} and  also \cite{Gautheron}. Moreover, 3-Lie algebras appeared in String Theory and
M-theory. In \cite{Basu},
 Basu and Harvey suggested to replace the Lie algebra appearing in the Nahm equation by a 3-Lie algebra for the  lifted Nahm equations.
Furthermore, in the context of Bagger-Lambert-Gustavsson model of multiple
M2-branes, Bagger-Lambert  managed to construct, using a ternary bracket,  an $N=2$
 supersymmetric version of the worldvolume theory of the M-theory membrane, see \cite{BL0} and also  \cite{BL3,BL2,HHM,P}.

Several algebraic aspects of $n$-Lie algebras  were studied in the last years. See   \cite{Realization,Baiclassification} for the construction, realization and classifications of 3-Lie algebras and $n$-Lie algebras.  Representation theory of $n$-Lie algebras was first introduced by Kasymov in \cite{Kasymov} and cohomologies were studied in \cite{cohomology}. The adjoint representation  is defined by the ternary bracket in which two elements are fixed.  Through fundamental objects one may also represent a 3-Lie algebra and more generally an $n$-Lie algebra by  a Leibniz algebra \cite{DT}. Following this approach,  deformations of 3-Lie algebras and $n$-Lie algebras are studied in   \cite{deformation,Tcohomology}, see  \cite{Makhlouf} for a review. In \cite{NR bracket of n-Lie}, the author defined a graded Lie algebra structure on the cochain complex of an $n$-Leibniz algebra and described an $n$-Leibniz structure as a canonical structure.  See  the review article \cite{review} for more details. In \cite{Liu-Makhlouf-Sheng}, the authors introduced the notion of a generalized representation of a 3-Lie algebra, by which abelian extensions of 3-Lie algebras are studied.

Due to its difficulty and less of tools, non-abelian extensions of 3-Lie algebras are not studied. In this paper, motivated by the work in \cite{nonabelin cohomology of Lie,Liu-Makhlouf-Sheng}, we find a suitable  approach which  uses  Maurer-Cartan elements to study non-abelian extensions of 3-Lie algebras.  We also show that the Leibniz algebra on the space of fundamental objects is a non-abelian extension of Leibniz algebras.

The paper is organized as follows. In Section 2, we give a review of non-abelian extensions of Leibniz algebras and cohomologies of   3-Lie algebras.  A characterization of    non-abelian extensions  of a 3-Lie algebra by another 3-Lie algebra  is given in Section 3 and several examples provided. In Section 4, we show that there is a one-to-one correspondence between isomorphism classes of non-abelian extensions of 3-Lie algebras and equivalence classes of Maurer-Cartan elements. Finally, we analyze in Section 5  the corresponding Leibniz algebra structure on the space of fundamental objects and show that it is a non-abelian extension of Leibniz algebras.

\section{Preliminaries}
In this paper, we work over an algebraically closed field $\K$ of characteristic 0 and all the vector spaces are over $\K$.

\subsection{Non-abelian extensions of Leibniz algebras}

A {\bf Leibniz algebra} is a vector space $\frkk$ endowed with a linear map $[\cdot,\cdot]_{\frkk}:\frkk\otimes\frkk\lon\frkk$ satisfying
\begin{eqnarray}
[x,[y,z]_{\frkk}]_{\frkk}=[[x,y]_{\frkk},z]_{\frkk}+[y,[x,z]_{\frkk}]_{\frkk},\,\,\,\,\forall x,y,z\in\frkk.
\end{eqnarray}
This is in fact a left Leibniz algebra. In this paper, we only consider left Leibniz algebras which we call Leibniz algebras. \vspace{3mm}

Let $(\frkk,[\cdot,\cdot]_\frkk)$ be a Leibniz algebra. We  denote by $\Der^L(\frkk)$ and $\Der^R(\frkk)$ the set of left derivations and the set of right derivations of $\g$ respectively:
\begin{eqnarray*}
\Der^L(\frkk)&=&\{D\in \gl(\frkk)\vert D[x,y]_{\frkk}=[Dx,y]_{\frkk}+[x,Dy]_{\frkk},\forall x,y\in \frkk\},\\
 \Der^R(\frkk)&=&\{D\in \gl(\frkk)\vert D[x,y]_{\frkk}=[x,Dy]_{\frkk}-[y,Dx]_{\frkk},\forall x,y\in \frkk\}.
 \end{eqnarray*}
Note that the right derivations are called anti-derivations in  \cite{Loday,Loday and Pirashvili}. It is easy to see that for all $x\in\frkk$, $\ad^L_x:\frkk\longrightarrow\frkk$, which is given by $\ad^L_x(y)=[x,y]_{\frkk}$, is a left derivation; $\ad^R_x:\frkk\longrightarrow\frkk$, which is given by $\ad^R_x(y)=[y,x]_{\frkk}$, is a right derivation.

\begin{defi}\label{defi:isomorphic}
\begin{itemize}
\item[\rm(1)] Let $\frkk$, $\frks$, $\hat{\frkk}$ be  Leibniz algebras. A non-abelian extension of Leibniz algebras is a short exact sequence of Leibniz algebras:
$$ 0\longrightarrow\frks\stackrel{i}{\longrightarrow}\hat{\frkk}\stackrel{p}\longrightarrow\frkk\longrightarrow0.$$
We say that $\hat{\frkk}$ is a non-abelian  extension of $\frkk$ by $\frks$.
\item[\rm(2)] A linear section of $\hat{\frkk}$ is a linear map $\sigma:\frkk\rightarrow\hat{\frkk}$ such that $p\circ \sigma=\Id$.
\end{itemize}
\end{defi}

Let $\hat{\frkk}$ be a non-abelian  extension of $\frkk$ by $\frks$, and  $\sigma:\frkk\rightarrow\widehat{\frkk}$ a linear section. Define $\omega:\frkk\otimes\frkk\rightarrow\frks$,
 $l:\frkk\longrightarrow\gl(\frks)$ and $r:\frkk\longrightarrow\gl(\frks)$ respectively by
\begin{eqnarray}
  \label{eq:str1}\omega(x,y)&=&[\sigma(x),\sigma(y)]_{\hat{\frkk}}-\sigma[x,y]_{\frkk},\quad \forall x,y\in\frkk,\\
 \label{eq:str2} l_x(\beta)&=&[\sigma(x),\beta]_{\hat{\frkk}},\quad \forall x\in\frkk,~\beta\in \frks,\\
 \label{eq:str31} r_y(\alpha)&=&[\alpha,\sigma(y)]_{\hat{\frkk}},\quad \forall y\in\frkk,~\alpha\in \frks.
\end{eqnarray}

Given a linear section, we have $\hat{\frkk}\cong\frkk\oplus \frks$ as vector spaces, and
the Leibniz algebra structure on  $\hat{\frkk}$ can be transferred to  $\frkk\oplus\frks$:
\begin{eqnarray}\label{bracket of nbe}
[x+\alpha,y+\beta]_{(l,r,\omega)}=[x,y]_{\frkk}+\omega(x,y)+l_x\beta+r_y\alpha+[\alpha,\beta]_{\frks}.\label{bracket of Ext}
\end{eqnarray}

\begin{pro}\label{non-abelian extension of LB}
With the above notations, $(\frkk\oplus\frks,[\cdot,\cdot]_{(l,r,\omega)})$ is a Leibniz algebra if and only if $l,r,\omega$ satisfy the following equalities:
\begin{eqnarray}
l_x[\alpha,\beta]_{\frks}&=&[l_x\alpha,\beta]_{\frks}+[\alpha,l_x\beta]_{\frks},\label{l Der}\\
r_x[\alpha,\beta]_{\frks}&=&[\alpha,r_x\beta]_{\frks}-[\beta,r_x\alpha]_{\frks},\label{r Der}\\
{[l_x\alpha+r_x\alpha,\beta]}_{\frks}&=&0,\label{Left Cen1}\\
{[l_x,l_y]}-l_{[x,y]_{\frkk}}&=&\ad^L_{\omega(x,y)},\label{rep1}\\
{[l_x,r_y]}-r_{[x,y]_{\frkk}}&=&\ad^R_{\omega(x,y)},\label{rep2}\\
r_y(r_x(\alpha)+l_x(\alpha))&=&0,\label{Left Cen2}\\
l_x\omega(y,z)-l_y\omega(x,z)-r_z\omega(x,y)&=&\omega([x,y]_{\frkk},z)-\omega(x,[y,z]_{\frkk})+\omega(y,[x,z]_{\frkk}).\label{cocycle}
\end{eqnarray}
\end{pro}

Eq. \eqref{l Der} means that $l_x\in\Der^L(\frks)$ and Eq. \eqref{r Der} means that $r_x\in\Der^R(\frks)$. See \cite{LSW} for more details about non-abelian extensions of Leibniz algebras.

\subsection{3-Lie algebras and their representations}
\begin{defi}
  A $3$-Lie algebra is a vector space $\g$ together with a skew-symmetric linear map  $[\cdot,\cdot,\cdot]_\g:
\wedge^3 \g\rightarrow \g$ such that the following {\bf Fundamental Identity (FI)} holds:
\begin{eqnarray}
\nonumber &&F_{x_1,x_2,x_3,x_4,x_5}\\
\nonumber&\triangleq&[x_1,x_2,[x_3,x_4,x_5]_\g]_\g-[[x_1,x_2,x_3]_\g,x_4,x_5]_\g-[x_3,[x_1,x_2,x_4]_\g,x_5]_\g-[x_3,x_4,[x_1,x_2,x_5]_\g]_\g\\
\label{eq:de1}&=&0.
\end{eqnarray}
\end{defi}

\begin{defi}
A  morphism of $3$-Lie algebras $f:(\g,[\cdot,\cdot,\cdot]_\g)\lon(\h,[\cdot,\cdot,\cdot]_\h)$ is a linear map $f:\g\lon\h$ such that
\begin{eqnarray}
f[x,y,z]_\g=[f(x),f(y),f(z)]_\h.
\end{eqnarray}
\end{defi}

Elements in $\wedge^2\g$ are called {\bf fundamental objects} of the $3$-Lie algebra $(\g,[\cdot,\cdot,\cdot]_\g)$. There is a bilinear operation $[\cdot,\cdot]_{\rm F}$ on $  \wedge^{2}\g$, which is given by
\begin{equation}\label{eq:bracketfunda}
~[\frkX,\frkY]_{\rm F}=[x_1,x_2,y_1]_\g\wedge y_2+y_1\wedge[x_1,x_2,y_2]_\g,\quad \forall \frkX=x_1\wedge x_2, ~\frkY=y_1\wedge y_2.
\end{equation}
It is well-known that $(\wedge^2\g,[\cdot,\cdot]_{\rm F})$ is a Leibniz algebra \cite{DT}, which plays  an important role in the theory of 3-Lie algebras.

\begin{defi}{\rm (\cite{Kasymov})}\label{defi:usualrep}
 A representation $\rho$ of a $3$-Lie algebra $\frkg$ on a vector space  $V$ is a linear map $\rho:\wedge^2\frkg\longrightarrow \gl( V),$ such that
\begin{eqnarray*}
  \rho(x_1,x_2)\rho(x_3,x_4)&=&\rho([x_1,x_2,x_3]_\g,x_4)+\rho(x_3,[x_1,x_2,x_4]_\g)+\rho(x_3,x_4)\rho(x_1,x_2),\\
\rho(x_1,[x_2,x_3,x_4]_\g)&=&\rho(x_3,x_4)\rho(x_1,x_2)-\rho(x_2,x_4)\rho(x_1,x_3)+\rho(x_2,x_3)\rho(x_1,x_4).
\end{eqnarray*}
\end{defi}

\begin{lem}\label{lem:semidirectp}
Let $\g$ be a $3$-Lie algebra, $V$  a vector space and $\rho:
\wedge^2\g\rightarrow \gl(V)$  a skew-symmetric linear
map. Then $(V;\rho)$ is a representation of $\g$ if and only if there
is a $3$-Lie algebra structure $($called the semidirect product$)$
on the direct sum of vector spaces  $\g\oplus V$, defined by
\begin{equation}\label{eq:sum}
[x_1+v_1,x_2+v_2,x_3+v_3]_{\rho}=[x_1,x_2,x_3]_\g+\rho(x_1,x_2)v_3+\rho(x_2,x_3)v_1+\rho(x_3,x_1)v_2,
\end{equation}
for $x_i\in \g, v_i\in V, 1\leq i\leq 3$. We denote this semidirect product $3$-Lie algebra by $\g\ltimes_\rho V.$
\end{lem}

A $p$-cochain on  $\frkg$ with  coefficients in a representation   $(V;\rho)$ is a linear map $$\alpha:\wedge^2\frkg\otimes\stackrel{(p-1)}{\cdots}\otimes\wedge^2\frkg\wedge\frkg\longrightarrow V.$$
Denote the space of $p$-cochains by $C^{p-1}(\g,V).$ The coboundary operator $\delta:C^{p-1}(\g,V)\longrightarrow C^{p}(\g,V)$  is given by
\begin{eqnarray}
\nonumber&&(\delta\alpha)(\frkX_1,\cdots ,\frkX_p,z)\\
\nonumber&=& \sum_{1\leq j<k\leq p}(-1)^j\alpha(\frkX_1,\cdots ,\hat{\frkX}_j,\cdots ,\frkX_{k-1},[\frkX_j,\frkX_k]_{\rm F},\frkX_{k+1},\cdots ,\frkX_{p},z)\\
\nonumber&&+\sum_{j=1}^p(-1)^j\alpha(\frkX_1,\cdots ,\hat{\frkX}_j,\cdots ,\frkX_{p},[\frkX_j,z]_\g)\\
\nonumber&&+\sum_{j=1}^p(-1)^{j+1}\rho(\frkX_j)\alpha(\frkX_1,\cdots ,\hat{\frkX}_j,\cdots ,\frkX_{p},z)\\
\label{eq:drho}&&+(-1)^{p+1}\Big(\rho(y_{p},z)\alpha(\frkX_1,\cdots ,\frkX_{p-1},x_{p} ) +\rho(z,x_{p})\alpha(\frkX_1,\cdots ,\frkX_{p-1},y_{p} ) \Big).
\end{eqnarray}

\section{Non-abelian extensions of $3$-Lie algebras}\label{sec:pre}
\begin{defi}
A non-abelian extension of a $3$-Lie algebra $(\g,[\cdot,\cdot,\cdot]_\g)$ by a $3$-Lie algebra $(\frkh,[\cdot,\cdot,\cdot]_{\frkh})$ is a short exact sequence of $3$-Lie algebra morphisms:
$$ 0\longrightarrow\h\stackrel{\iota}{\longrightarrow}\hat{\g}\stackrel{p}\longrightarrow\g\longrightarrow0,$$
where $(\hat{\mathfrak{g}},[\cdot,\cdot,\cdot]_{\hat{\g}})$ is a $3$-Lie algebra.
\end{defi}

\begin{defi}\label{defi:iso}
Two extensions of $\mathfrak{g}$ by $\mathfrak{h}$, $(\hat{\g}_1,[\cdot,\cdot,\cdot]_{\hat{\g}_1})$ and $(\hat{\g}_2,[\cdot,\cdot,\cdot]_{\hat{\g}_2})$, are said to be isomorphic if there exists a $3$-Lie algebra morphism $\theta:\hat{\mathfrak{g}}_{2}\lon \hat{\mathfrak{g}}_{1}$ such that we have the following commutative diagram:
\label{iso}\[\begin{CD}
0@>>>\mathfrak{h}@>\iota_{2}>>\hat{\mathfrak{g}}_{2}@>p_{2}>>{\mathfrak{g}}             @>>>0\\
@.    @|                       @V\theta VV                     @|                       @.\\
0@>>>\mathfrak{h}@>\iota_{1}>>\hat{\mathfrak{g}}_{1}@>p_{1}>>{\mathfrak{g}}             @>>>0
.\end{CD}\]
\end{defi}
Given a section $\sigma$ of $\hat{\g}$, define $\rho:\wedge^2\g\longrightarrow\gl(\h)$, $\nu:\g\longrightarrow\Hom(\wedge^2\h,\h)$ and $\omega:\wedge^3\frkg\longrightarrow\frkh$  by
\begin{eqnarray*}
 \rho(x,y)(u)&=&[\sigma(x),\sigma(y),u]_{\hat{\frkg}},\\
  \nu(x)(u,v)&=&[\sigma(x),u,v]_{\hat{\frkg}},\\
\omega(x,y,z)&=&[\sigma(x),\sigma(y),\sigma(z)]_{\hat{\frkg}}-\sigma[x,y,z]_{\frkg}.
\end{eqnarray*}
Obviously, $\hat{\g}$ is isomorphic to $\g\oplus\h$ as vector spaces. Transfer the $3$-Lie algebra structure on $\hat{\mathfrak{g}}$ to that on $\mathfrak{g}\oplus \mathfrak{h}$, we obtain a $3$-Lie algebra $(\mathfrak{g}\oplus \mathfrak{h},[\cdot,\cdot,\cdot]_{(\rho,\nu,\omega)})$, where $[\cdot,\cdot]_{(\rho,\nu,\omega)}$ is given by
\begin{eqnarray*}
&&[x_1+v_1,x_2+v_2,x_3+v_3]_{(\rho,\nu,\omega)}\\&&=[x_1,x_2,x_3]_{\frkg}+\omega(x_1,x_2,x_3)+\rho(x_1,x_2)v_3+\rho(x_2,x_3)v_1+\rho(x_3,x_1)v_2\\
&&+\nu(x_1)(v_2,v_3)+\nu(x_2)(v_3,v_1)+\nu(x_3)(v_1,v_2)+[v_1,v_2,v_3]_{\frkh}.
\end{eqnarray*}

The following proposition gives the conditions on $\rho,\nu$ and $\omega$ such that $(\g\oplus\h, [\cdot,\cdot,\cdot]_{(\rho,\nu,\omega)})$ is a $3$-Lie algebra.
\begin{pro}\label{3-lie}
With the above notations, $(\g\oplus\h, [\cdot,\cdot,\cdot]_{(\rho,\nu,\omega)})$ is a $3$-Lie algebra if and only if  $\rho,\nu$ and $\omega $ satisfy, for all $x_1,\cdots,x_5\in \g$, $ v_1\cdots, v_5\in\h$, the following conditions
\begin{eqnarray}
\nonumber0&=&\omega(x_1,x_2,[x_3,x_4,x_5]_\g)+\rho(x_1,x_2)\omega(x_3,x_4,x_5)-\omega([x_1,x_2,x_3]_\g,x_4,x_5)\\
\nonumber&&-\rho(x_4,x_5)\omega(x_1,x_2,x_3)-\omega(x_3,[x_1,x_2,x_4]_\g,x_5)+\rho(x_3,x_5)\omega(x_1,x_2,x_4)\\
\label{p1}&&-\omega(x_3,x_4,[x_1,x_2,x_5]_\g)-\rho(x_3,x_4)\omega(x_1,x_2,x_5),\\
\nonumber0&=&\rho(x_1,x_2)(\rho(x_3,x_4)v_5)-\rho([x_1,x_2,x_3]_\g,x_4)v_5+\nu(x_4)(\omega(x_1,x_2,x_3),v_5)\\
\label{p2}&&-\rho(x_3,[x_1,x_2,x_4]_\g)v_5-\nu(x_3)(\omega(x_1,x_2,x_4),v_5)-\rho(x_3,x_4)(\rho(x_1,x_2)v_5),\\
\nonumber0&=&-\rho(x_1,[x_3,x_4,x_5]_\g)v_2+\nu(x_1)(v_2,\omega(x_3,x_4,x_5))+\rho(x_4,x_5)(\rho(x_1,x_3)v_2)\\
\label{p3}&&-\rho(x_3,x_5)(\rho(x_1,x_4)v_2)+\rho(x_3,x_4)(\rho(x_1,x_5)v_2),\\
\nonumber0&=&\rho(x_1,x_2)(\nu(x_3)(v_4,v_5))-\nu([x_1,x_2,x_3]_\g)(v_4,v_5)-[\omega(x_1,x_2,x_3),v_4,v_5]_\h\\
\label{p4}&&-\nu(x_3)(\rho(x_1,x_2)v_4,v_5)-\nu(x_3)(v_4,\rho(x_1,x_2)v_5),\\
\nonumber0&=&\nu(x_1)(v_2,\rho(x_3,x_4)v_5)-\nu(x_4)(\rho(x_1,x_3)v_2,v_5)+\nu(x_3)(\rho(x_1,x_4)v_2,v_5)\\
\label{p5}&&-\rho(x_3,x_4)(\nu(x_1)(v_2,v_5)),\\
\nonumber0&=&\nu([x_3,x_4,x_5]_\g)(v_1,v_2)+[v_1,v_2,\omega(x_3,x_4,x_5)]_\h-\rho(x_4,x_5)(\nu(x_3)(v_1,v_2))\\
\label{p6}&&+\rho(x_3,x_5)(\nu(x_4)(v_1,v_2))-\rho(x_3,x_4)(\nu(x_5)(v_1,v_2)),\\
\nonumber0&=&\rho(x_1,x_2)[v_3,v_4,v_5]_\h-[\rho(x_1,x_2)v_3,v_4,v_5]_\h-[v_3,\rho(x_1,x_2)v_4,v_5]_\h\\
\label{p7}&&-[v_3,v_4,\rho(x_1,x_2)v_5]_\h,\\
\nonumber0&=&\nu(x_1)(v_2,\nu(x_3)(v_4,v_5))+[\rho(x_1,x_3)v_2,v_4,v_5]_\h-\nu(x_3)(\nu(x_1)(v_2,v_4),v_5)\\
\label{p8}&&-\nu(x_3)(v_4,\nu(x_1)(v_2,v_5)),\\
\nonumber0&=&[v_1,v_2,\rho(x_3,x_4)v_5]_\h+\nu(x_4)(\nu(x_3)(v_1,v_2),v_5)-\nu(x_3)(\nu(x_4)(v_1,v_2),v_5)\\
\label{p9}&&-\rho(x_3,x_4)[v_1,v_2,v_5]_\h,\\
\nonumber0&=&\nu(x_1)(v_2,[v_3,v_4,v_5]_\h)-[\nu(x_1)(v_2,v_3),v_4,v_5]_\h-[v_3,\nu(x_1)(v_2,v_4),v_5]_\h\\
\label{p10}&&-[v_3,v_4,\nu(x_1)(v_2,v_5)]_\h,\\
\nonumber0&=&[v_1,v_2,\nu(x_3)(v_4,v_5)]_\h-[\nu(x_3)(v_1,v_2),v_4,v_5]_\h-\nu(x_3)([v_1,v_2,v_4]_\h,v_5)\\
\label{p11}&&-\nu(x_3)(v_4,[v_1,v_2,v_5]_\h).
\end{eqnarray}
\end{pro}

\pf Assume that $(\g\oplus\h, [\cdot,\cdot,\cdot]_{(\rho,\nu,\omega)})$ is a $3$-Lie algebra. By
\begin{eqnarray*}
[x_1,x_2,[x_3,x_4,x_5]_{(\rho,\nu,\omega)}]_{(\rho,\nu,\omega)}&=&[[x_1,x_2,x_3]_{(\rho,\nu,\omega)},x_4,x_5]_{(\rho,\nu,\omega)}
+[x_3,[x_1,x_2,x_4]_{(\rho,\nu,\omega)},x_5]_{(\rho,\nu,\omega)}\\
&&+[x_3,x_4,[x_1,x_2,x_5]_{(\rho,\nu,\omega)}]_{(\rho,\nu,\omega)},
\end{eqnarray*}
we deduce that \eqref{p1} holds. By
\begin{eqnarray*}
[x_1,x_2,[x_3,x_4,v_5]_{(\rho,\nu,\omega)}]_{(\rho,\nu,\omega)}&=&[[x_1,x_2,x_3]_{(\rho,\nu,\omega)},x_4,v_5]_{(\rho,\nu,\omega)}
+[x_3,[x_1,x_2,x_4]_{(\rho,\nu,\omega)},v_5]_{(\rho,\nu,\omega)}\\
&&+[x_3,x_4,[x_1,x_2,v_5]_{(\rho,\nu,\omega)}]_{(\rho,\nu,\omega)},
\end{eqnarray*}
we deduce that \eqref{p2} holds. By
\begin{eqnarray*}
[x_1,v_2,[x_3,x_4,x_5]_{(\rho,\nu,\omega)}]_{(\rho,\nu,\omega)}&=&[[x_1,v_2,x_3]_{(\rho,\nu,\omega)},x_4,x_5]_{(\rho,\nu,\omega)}
+[x_3,[x_1,v_2,x_4]_{(\rho,\nu,\omega)},x_5]_{(\rho,\nu,\omega)}\\
&&+[x_3,x_4,[x_1,v_2,x_5]_{(\rho,\nu,\omega)}]_{(\rho,\nu,\omega)},
\end{eqnarray*}
we deduce that \eqref{p3} holds. By
\begin{eqnarray*}
[x_1,x_2,[x_3,v_4,v_5]_{(\rho,\nu,\omega)}]_{(\rho,\nu,\omega)}&=&[[x_1,x_2,x_3]_{(\rho,\nu,\omega)},v_4,v_5]_{(\rho,\nu,\omega)}
+[x_3,[x_1,x_2,v_4]_{(\rho,\nu,\omega)},v_5]_{(\rho,\nu,\omega)}\\
&&+[x_3,v_4,[x_1,x_2,v_5]_{(\rho,\nu,\omega)}]_{(\rho,\nu,\omega)},
\end{eqnarray*}
we deduce that \eqref{p4} holds. By
\begin{eqnarray*}
[x_1,v_2,[x_3,x_4,v_5]_{(\rho,\nu,\omega)}]_{(\rho,\nu,\omega)}&=&[[x_1,v_2,x_3]_{(\rho,\nu,\omega)},x_4,v_5]_{(\rho,\nu,\omega)}
+[x_3,[x_1,v_2,x_4]_{(\rho,\nu,\omega)},v_5]_{(\rho,\nu,\omega)}\\
&&+[x_3,x_4,[x_1,v_2,v_5]_{(\rho,\nu,\omega)}]_{(\rho,\nu,\omega)},
\end{eqnarray*}
we deduce that \eqref{p5} holds. By
\begin{eqnarray*}
[v_1,v_2,[x_3,x_4,x_5]_{(\rho,\nu,\omega)}]_{(\rho,\nu,\omega)}&=&[[v_1,v_2,x_3]_{(\rho,\nu,\omega)},x_4,x_5]_{(\rho,\nu,\omega)}
+[x_3,[v_1,v_2,x_4]_{(\rho,\nu,\omega)},x_5]_{(\rho,\nu,\omega)}\\
&&+[x_3,x_4,[v_1,v_2,x_5]_{(\rho,\nu,\omega)}]_{(\rho,\nu,\omega)},
\end{eqnarray*}
we deduce that \eqref{p6} holds. By
\begin{eqnarray*}
[x_1,x_2,[v_3,v_4,v_5]_{(\rho,\nu,\omega)}]_{(\rho,\nu,\omega)}&=&[[x_1,x_2,v_3]_{(\rho,\nu,\omega)},v_4,v_5]_{(\rho,\nu,\omega)}
+[v_3,[x_1,x_2,v_4]_{(\rho,\nu,\omega)},v_5]_{(\rho,\nu,\omega)}\\
&&+[v_3,v_4,[x_1,x_2,v_5]_{(\rho,\nu,\omega)}]_{(\rho,\nu,\omega)},
\end{eqnarray*}
we deduce that \eqref{p7} holds. By
\begin{eqnarray*}
[x_1,v_2,[x_3,v_4,v_5]_{(\rho,\nu,\omega)}]_{(\rho,\nu,\omega)}&=&[[x_1,v_2,x_3]_{(\rho,\nu,\omega)},v_4,v_5]_{(\rho,\nu,\omega)}
+[x_3,[x_1,v_2,v_4]_{(\rho,\nu,\omega)},v_5]_{(\rho,\nu,\omega)}\\
&&+[x_3,v_4,[x_1,v_2,v_5]_{(\rho,\nu,\omega)}]_{(\rho,\nu,\omega)},
\end{eqnarray*}
we deduce that \eqref{p8} holds. By
\begin{eqnarray*}
[v_1,v_2,[x_3,x_4,v_5]_{(\rho,\nu,\omega)}]_{(\rho,\nu,\omega)}&=&[[v_1,v_2,x_3]_{(\rho,\nu,\omega)},x_4,v_5]_{(\rho,\nu,\omega)}
+[x_3,[v_1,v_2,x_4]_{(\rho,\nu,\omega)},v_5]_{(\rho,\nu,\omega)}\\
&&+[x_3,x_4,[v_1,v_2,v_5]_{(\rho,\nu,\omega)}]_{(\rho,\nu,\omega)},
\end{eqnarray*}
we deduce that \eqref{p9} holds. By
\begin{eqnarray*}
[x_1,v_2,[v_3,x_4,v_5]_{(\rho,\nu,\omega)}]_{(\rho,\nu,\omega)}&=&[[x_1,v_2,v_3]_{(\rho,\nu,\omega)},v_4,v_5]_{(\rho,\nu,\omega)}
+[v_3,[x_1,v_2,v_4]_{(\rho,\nu,\omega)},v_5]_{(\rho,\nu,\omega)}\\
&&+[v_3,v_4,[x_1,v_2,v_5]_{(\rho,\nu,\omega)}]_{(\rho,\nu,\omega)},
\end{eqnarray*}
we deduce that \eqref{p10} holds. By
\begin{eqnarray*}
[v_1,v_2,[x_3,v_4,v_5]_{(\rho,\nu,\omega)}]_{(\rho,\nu,\omega)}&=&[[v_1,v_2,x_3]_{(\rho,\nu,\omega)},v_4,v_5]_{(\rho,\nu,\omega)}
+[x_3,[v_1,v_2,v_4]_{(\rho,\nu,\omega)},v_5]_{(\rho,\nu,\omega)}\\
&&+[x_3,v_4,[v_1,v_2,v_5]_{(\rho,\nu,\omega)}]_{(\rho,\nu,\omega)},
\end{eqnarray*}
we deduce that \eqref{p11} holds.

Conversely, if \eqref{p1}-\eqref{p11} hold, it is straightforward to see that $(\g\oplus\h, [\cdot,\cdot,\cdot]_{(\rho,\nu,\omega)})$ is a $3$-Lie algebra. \qed

\begin{ex}
Let $\g$ be the simple  $4$-dimensional $3$-Lie algebra defined with respect to a  basis $\{x_1,x_2,x_3,x_4\}$ by the skew-symmetric brackets
\begin{eqnarray*}
[x_1,x_2,x_3]=x_4, \; [x_1,x_2,x_4]=x_3, \; [x_1,x_3,x_4]=x_2, \; [x_2,x_3,x_4]=x_1,
\end{eqnarray*} and let $\mathfrak{h}$ be the $3$-dimensional $ 3$-Lie algebra defined with respect   to basis $\{v_1,v_2,v_3\}$ by  $$[v_1,v_2,v_3]=v_1.$$ Then every non-abelian  extension of $\g$ by $\mathfrak{h}$ is given by $\rho =0$.  The following families of  $\nu$ and $\omega$ provide   non-abelian  extensions of $\g$ by $\mathfrak{h}$

\begin{itemize}
\item $ \nu(x_1)(v_1,v_2)=r_1v_1,  \; \nu(x_1)(v_1,v_3)=r_2v_1,  \; \nu(x_1)(v_2,v_3)=r_3v_1,  \;$\\ $
 \nu(x_2)(v_1,v_2)=\frac{r_1r_4}{r_2} v_1,  \; \nu(x_2)(v_1,v_3)=r_4v_1,  \;
 \nu(x_2)(v_2,v_3)=\frac{r_3r_4}{r_2}v_1,  \;$\\ $
 \nu(x_3)(v_1,v_2)=\frac{r_1r_5}{r_2}v_1,  \; \nu(x_3)(v_1,v_3)=r_5v_1,  \;
 \nu(x_3)(v_2,v_3)=\frac{r_3r_5}{r_2}v_1, \;
 $\\ $
 \nu(x_4)(v_1,v_2)=\frac{r_1r_6}{r_2}v_1,  \; \nu(x_4)(v_1,v_3)=r_6v_1,  \;
 \nu(x_4)(v_2,v_3)=\frac{r_3r_6}{r_2}v_1, $

\item $\omega (x_1,x_2,x_3)=-\frac{r_3r_6}{r_2}v_1+r_6v_2-\frac{r_1r_6}{r_2}v_3,$\\
$\omega (x_1,x_2,x_4)=-\frac{r_3r_5}{r_2}v_1+r_5v_2-\frac{r_1r_5}{r_2}v_3,$\\
$\omega (x_1,x_3,x_4)=-\frac{r_3r_4}{r_2}v_1+r_4v_2-\frac{r_1r_4}{r_2}v_3,$\\
$\omega (x_2,x_3,x_4)=-r_3v_1+r_2v_2-r_1v_3,$
\end{itemize}

where $r_i$  are parameters in $\K$.

\end{ex}

\begin{ex}
Now, let $\g$ be the $3$-dimensional $3$-Lie algebra defined with respect to a  basis $\{x_1,x_2,x_3\}$ by the skew-symmetric bracket $[x_1,x_2,x_3]=x_1$ and let $\mathfrak{h}$ be the same $3$-Lie algebra which we consider with respect to basis $\{v_1,v_2,v_3\}$, that is $[v_1,v_2,v_3]=v_1$. Then every non-abelian  extension of $\g$ by $\mathfrak{h}$ is given by one of the following triples $(\rho,\nu,\omega)$.

\begin{enumerate}
\item
\begin{itemize}
\item $   \rho(x_1,x_2)(v_1)=0,  \;
\rho(x_1,x_2)(v_2)=0 ,   \;
\rho(x_1,x_2)(v_3)=0 ,   \;$\\ $
 \rho(x_1,x_3)(v_1)=0, \;
  \rho(x_1,x_3)(v_2)=0, \;
    \rho(x_1,x_3)(v_3)=0, \;$\\ $
 \rho(x_2,x_3)(v_1)= (r_3r_6-r_4r_5)v_1,    \;  \rho(x_2,x_3)(v_2)=r_1 v_1,\;  \rho(x_2,x_3)(v_3)=r_2 v_1,$

\item $ \nu(x_1)(v_1,v_2)=0,  \; \nu(x_1)(v_1,v_3)=0,  \; \nu(x_1)(v_2,v_3)=0,  \;$\\ $
 \nu(x_2)(v_1,v_2)=r_3 v_1,  \; \nu(x_2)(v_1,v_3)=r_4v_1,  \;
 \nu(x_2)(v_2,v_3)=\frac{r_2r_3-r_1r_4}{r_4r_5-r_3r_6}v_1,  \;$\\ $
 \nu(x_3)(v_1,v_2)=r_5v_1,  \; \nu(x_3)(v_1,v_3)=r_6v_1,  \;
 \nu(x_3)(v_2,v_3)=\frac{r_2r_5-r_1r_6}{r_4r_5-r_3r_6}v_1, $

\item $\omega (x_1,x_2,x_3)=0,$
\end{itemize}
\item
\begin{itemize}
\item $   \rho(x_1,x_2)(v_1)=0,  \;
\rho(x_1,x_2)(v_2)=0 ,   \;
\rho(x_1,x_2)(v_3)=0 ,   \;$\\ $
 \rho(x_1,x_3)(v_1)=0, \;
  \rho(x_1,x_3)(v_2)=0, \;
    \rho(x_1,x_3)(v_3)=0, \;$\\ $
 \rho(x_2,x_3)(v_1)= -r_3r_4 v_1,    \;  \rho(x_2,x_3)(v_2)=r_1 v_1,\;  \rho(x_2,x_3)(v_3)=r_2 v_1,$

\item $ \nu(x_1)(v_1,v_2)=0,  \; \nu(x_1)(v_1,v_3)=0,  \; \nu(x_1)(v_2,v_3)=0,  \;$\\ $
 \nu(x_2)(v_1,v_2)=0,  \; \nu(x_2)(v_1,v_3)=r_3v_1,  \;
 \nu(x_2)(v_2,v_3)=\frac{r_1}{r_4}v_1,  \;$\\ $
 \nu(x_3)(v_1,v_2)=r_4v_1,  \; \nu(x_3)(v_1,v_3)=r_5v_1,  \;
 \nu(x_3)(v_2,v_3)=\frac{r_2r_4-r_1r_5}{r_3r_4}v_1, $

\item $\omega (x_1,x_2,x_3)=0,$
\end{itemize}
\item
\begin{itemize}
\item $   \rho(x_1,x_2)(v_1)=0,  \;
\rho(x_1,x_2)(v_2)=0 ,   \;
\rho(x_1,x_2)(v_3)=0 ,   \;$\\ $
 \rho(x_1,x_3)(v_1)=0, \;
  \rho(x_1,x_3)(v_2)=0, \;
    \rho(x_1,x_3)(v_3)=0, \;$\\ $
 \rho(x_2,x_3)(v_1)=0,    \;  \rho(x_2,x_3)(v_2)=r_1 v_1,\;  \rho(x_2,x_3)(v_3)=\frac{r_3r_1}{r_2} v_1,$

\item $ \nu(x_1)(v_1,v_2)=0,  \; \nu(x_1)(v_1,v_3)=0,  \; \nu(x_1)(v_2,v_3)=0,  \;$\\ $
 \nu(x_2)(v_1,v_2)=r_2v_1,  \; \nu(x_2)(v_1,v_3)=r_3v_1,  \;
 \nu(x_2)(v_2,v_3)=r_4v_1,  \;$\\ $
 \nu(x_3)(v_1,v_2)=r_5v_1,  \; \nu(x_3)(v_1,v_3)=\frac{r_3r_5}{r_2}v_1,  \;
 \nu(x_3)(v_2,v_3)=\frac{r_4r_5+r_1}{r_2}v_1, $

\item $\omega (x_1,x_2,x_3)=0,$
\end{itemize}
\item
\begin{itemize}
\item $   \rho=0$

\item $ \nu(x_1)(v_1,v_2)=r_1v_1,  \; \nu(x_1)(v_1,v_3)=r_2v_1,  \; \nu(x_1)(v_2,v_3)=r_3v_1,  \;$\\ $
 \nu(x_2)(v_1,v_2)=r_4v_1,  \; \nu(x_2)(v_1,v_3)=\frac{r_2r_4}{r_1}v_1,  \;
 \nu(x_2)(v_2,v_3)=\frac{r_3r_4}{r_1}v_1,  \;$\\ $
 \nu(x_3)(v_1,v_2)=r_5v_1,  \; \nu(x_3)(v_1,v_3)=\frac{r_2r_5}{r_1}v_1,  \;
 \nu(x_3)(v_2,v_3)=\frac{r_3r_5}{r_1}v_1, $

\item $\omega (x_1,x_2,x_3)=-r_3 v_1+r_2v_2-r_1v_3,$
\end{itemize}
\item
\begin{itemize}
\item $   \rho(x_1,x_2)(v_1)=0,  \;
\rho(x_1,x_2)(v_2)=0 ,   \;
\rho(x_1,x_2)(v_3)=0 ,   \;$\\ $
 \rho(x_1,x_3)(v_1)=0, \;
  \rho(x_1,x_3)(v_2)=0, \;
    \rho(x_1,x_3)(v_3)=0, \;$\\ $
 \rho(x_2,x_3)(v_1)=0,    \;  \rho(x_2,x_3)(v_2)=r_1 v_1,\;  \rho(x_2,x_3)(v_3)=\frac{r_3r_1}{r_2} v_1,$

\item $ \nu(x_1)(v_1,v_2)=0,  \; \nu(x_1)(v_1,v_3)=0,  \; \nu(x_1)(v_2,v_3)=0,  \;$\\ $
 \nu(x_2)(v_1,v_2)=0,  \; \nu(x_2)(v_1,v_3)=0,  \;
 \nu(x_2)(v_2,v_3)=-\frac{r_1}{r_2}v_1,  \;$\\ $
 \nu(x_3)(v_1,v_2)=r_2v_1,  \; \nu(x_3)(v_1,v_3)=r_3v_1,  \;
 \nu(x_3)(v_2,v_3)=r_4v_1, $

\item $\omega (x_1,x_2,x_3)=0,$
\end{itemize}
\item
\begin{itemize}
\item $   \rho(x_1,x_2)(v_1)=0,  \;
\rho(x_1,x_2)(v_2)=0 ,   \;
\rho(x_1,x_2)(v_3)=0 ,   \;$\\ $
 \rho(x_1,x_3)(v_1)=0, \;
  \rho(x_1,x_3)(v_2)=0, \;
    \rho(x_1,x_3)(v_3)=0, \;$\\ $
 \rho(x_2,x_3)(v_1)=0,    \;  \rho(x_2,x_3)(v_2)=0,\;  \rho(x_2,x_3)(v_3)=r_1 v_1,$

\item $ \nu(x_1)(v_1,v_2)=0,  \; \nu(x_1)(v_1,v_3)=0,  \; \nu(x_1)(v_2,v_3)=0,  \;$\\ $
 \nu(x_2)(v_1,v_2)=0,  \; \nu(x_2)(v_1,v_3)=r_2 v_1,  \;
 \nu(x_2)(v_2,v_3)=r_3v_1,  \;$\\ $
 \nu(x_3)(v_1,v_2)=0,  \; \nu(x_3)(v_1,v_3)=r_4v_1,  \;
 \nu(x_3)(v_2,v_3)=\frac{r_3r_4+r_1}{r_2} v_1, $

\item $\omega (x_1,x_2,x_3)=0,$
\end{itemize}

\item
\begin{itemize}
\item $   \rho=0$

\item $ \nu(x_1)(v_1,v_2)=0,  \; \nu(x_1)(v_1,v_3)=r_1v_1,  \; \nu(x_1)(v_2,v_3)=r_2v_1,  \;$\\ $
 \nu(x_2)(v_1,v_2)=0,  \; \nu(x_2)(v_1,v_3)=r_3v_1,  \;
 \nu(x_2)(v_2,v_3)=\frac{r_2r_3}{r_1}v_1,  \;$\\ $
 \nu(x_3)(v_1,v_2)=0,  \; \nu(x_3)(v_1,v_3)=r_4v_1,  \;
 \nu(x_3)(v_2,v_3)=\frac{r_2r_4}{r_1}v_1, $

\item $\omega (x_1,x_2,x_3)=-r_2 v_1+r_1v_2,$
\end{itemize}

\item
\begin{itemize}
\item $   \rho=0$

\item $ \nu(x_1)(v_1,v_2)=r_1 v_1,  \; \nu(x_1)(v_1,v_3)=r_2v_1,  \; \nu(x_1)(v_2,v_3)=0,  \;$\\ $
 \nu(x_2)(v_1,v_2)=r_3v_1,  \; \nu(x_2)(v_1,v_3)=\frac{r_2r_3}{r_1}v_1,  \;
 \nu(x_2)(v_2,v_3)=0,  \;$\\ $
 \nu(x_3)(v_1,v_2)=r_4 v_1,  \; \nu(x_3)(v_1,v_3)=\frac{r_2r_4}{r_1}v_1,  \;
 \nu(x_3)(v_2,v_3)=0, $

\item $\omega (x_1,x_2,x_3)=r_2 v_2-r_1v_3,$
\end{itemize}

\item
\begin{itemize}
\item $   \rho=0$

\item $ \nu(x_1)(v_1,v_2)=0,  \; \nu(x_1)(v_1,v_3)=0,  \; \nu(x_1)(v_2,v_3)=r_1v_1,  \;$\\ $
 \nu(x_2)(v_1,v_2)=0,  \; \nu(x_2)(v_1,v_3)=0,  \;
 \nu(x_2)(v_2,v_3)=r_2v_1,  \;$\\ $
 \nu(x_3)(v_1,v_2)=0,  \; \nu(x_3)(v_1,v_3)=0,  \;
 \nu(x_3)(v_2,v_3)=r_3v_1, $

\item $\omega (x_1,x_2,x_3)=-r_1v_1,$
\end{itemize}
\item
\begin{itemize}
\item $   \rho(x_1,x_2)(v_1)=0,  \;
\rho(x_1,x_2)(v_2)=0 ,   \;
\rho(x_1,x_2)(v_3)=0 ,   \;$\\ $
 \rho(x_1,x_3)(v_1)=0, \;
  \rho(x_1,x_3)(v_2)=0, \;
    \rho(x_1,x_3)(v_3)=0, \;$\\ $
 \rho(x_2,x_3)(v_1)=0,    \;  \rho(x_2,x_3)(v_2)=0,\;  \rho(x_2,x_3)(v_3)=r_1 v_1,$

\item $ \nu(x_1)(v_1,v_2)=0,  \; \nu(x_1)(v_1,v_3)=0,  \; \nu(x_1)(v_2,v_3)=0,  \;$\\ $
 \nu(x_2)(v_1,v_2)=0,  \; \nu(x_2)(v_1,v_3)=0,  \;
 \nu(x_2)(v_2,v_3)=-\frac{r_1}{r_2}v_1,  \;$\\ $
 \nu(x_3)(v_1,v_2)=0,  \; \nu(x_3)(v_1,v_3)=r_2v_1,  \;
 \nu(x_3)(v_2,v_3)=r_3v_1, $

\item $\omega (x_1,x_2,x_3)=0,$
\end{itemize}
\end{enumerate}
where $r_i$  are parameters in $\K$.

\end{ex}

\emptycomment{
\begin{rmk}
 \begin{itemize}
  \item[] Eq. \eqref{eq:n5} means that $\rho(x_1,x_2)$ is a derivation on $\frkh$, that is $\rho(x_1,x_2)\in \Der(\frkh),$ for all $x_1,x_2\in\h.$

  \item[] Eq. \eqref{eq:n6} means that $\nu(x)(v,\cdot)$ is a derivation on $\frkh$ for all $x\in\g,~v\in\h$.

  \item[] Eq. \eqref{eq:n7} means that $$
\nu(x)(\huaX\cdot \huaY)=-[\nu(x)(\huaX),\huaY]+[\huaX,\nu(x)(\huaY)],
$$
i.e. $\nu(x)$ is a derivation on the Leibniz algebra $\wedge^2\g\oplus \g$?.

\item[] Eq. \eqref{eq:n3} can be understood as
$$
[\nu(x)(u,\cdot),\nu(y)]+\ad_{\rho(x,y)(u)}=0.
$$

\item[] Eq. \eqref{eq:n4} can be understood as
$$
[\ad_{v_1,v_2},\rho(x_1,x_2)]=-\nu(x_2)(\nu(x_1)(v_1,v_2),\cdot)+\nu(x_1)(\nu(x_2)(v_1,v_2),\cdot).
$$
 \end{itemize}
\end{rmk}
}

Any non-abelian extension, by choosing a section,  is isomorphic to
$(\mathfrak{g}\oplus \mathfrak{h},[\cdot,\cdot]_{(\rho,\nu,\omega)})$. Therefore, we only consider in the sequel non-abelian extensions of the form $(\mathfrak{g}\oplus \mathfrak{h},[\cdot,\cdot]_{(\rho,\nu,\omega)})$.

\emptycomment{
\begin{pro}
Let $(\mathfrak{g}\oplus \mathfrak{h},[\cdot,\cdot,\cdot]_{(\rho^{1},\nu^1,\omega^{1})})$ and $(\mathfrak{g}\oplus \mathfrak{h},[\cdot,\cdot,\cdot]_{(\rho^{2},\nu^2,\omega^{2})})$ be  two non-abelian extensions of $\g$
by $\frkh$. Then the two extensions are isomorphic if and only if there is a linear map $\phi:\g\lon\frkh$ such that the following equalities holds:
\begin{eqnarray}
\label{iso1}\nu^2(x_1)(v_2,v_3)-\nu^1(x_1)(v_2,v_3)&=&[\phi(x_1),v_2,v_3]_\h,\\
\nonumber\rho^2(x_1,x_2)v_3-\rho^1(x_1,x_2)v_3&=&\nu^1(x_1)(\phi(x_2),v_3)+\nu^1(x_2)(v_3,\phi(x_1))\\
\label{iso2}&&+[\phi(x_1),\phi(x_2),v_3]_\h,\\
\nonumber\omega^2(x_1,x_2,x_3)-\omega^1(x_1,x_2,x_3)&=&\rho^1(x_1,x_2)\phi(x_3)+\rho^1(x_2,x_3)\phi(x_1)\\
\nonumber&&+\rho^1(x_3,x_1)\phi(x_2)+\nu^1(x_1)(\phi(x_2),\phi(x_3))\\
\nonumber&&+\nu^1(x_2)(\phi(x_3),\phi(x_1))+\nu^1(x_3)(\phi(x_1),\phi(x_2))\\
\label{iso3}&&+[\phi(x_1),\phi(x_2),\phi(x_3)]_\h-\phi[x_1,x_2,x_3]_\g.
\end{eqnarray}
\end{pro}
}

\begin{pro}
Let $(\mathfrak{g}\oplus \mathfrak{h},[\cdot,\cdot,\cdot]_{(\rho^{1},\nu^1,\omega^{1})})$ and $(\mathfrak{g}\oplus \mathfrak{h},[\cdot,\cdot,\cdot]_{(\rho^{2},\nu^2,\omega^{2})})$ be  two non-abelian extensions of $\g$
by $\frkh$. Then the two extensions are isomorphic if and only if there is a linear map $\xi:\g\lon\frkh$ such that the following equalities holds:
\begin{eqnarray}
\label{iso1}\nu^2(x_1)(v_2,v_3)-\nu^1(x_1)(v_2,v_3)&=&-[\xi(x_1),v_2,v_3]_\h,\\
\nonumber\rho^2(x_1,x_2)v_3-\rho^1(x_1,x_2)v_3&=&-\nu^1(x_1)(\xi(x_2),v_3)-\nu^1(x_2)(v_3,\xi(x_1))\\
\label{iso2}&&+[\xi(x_1),\xi(x_2),v_3]_\h,\\
\nonumber\omega^2(x_1,x_2,x_3)-\omega^1(x_1,x_2,x_3)&=&-\rho^1(x_1,x_2)\xi(x_3)-\rho^1(x_2,x_3)\xi(x_1)\\
\nonumber&&-\rho^1(x_3,x_1)\xi(x_2)+\nu^1(x_1)(\xi(x_2),\xi(x_3))\\
\nonumber&&+\nu^1(x_2)(\xi(x_3),\xi(x_1))+\nu^1(x_3)(\xi(x_1),\xi(x_2))\\
\label{iso3}&&-[\xi(x_1),\xi(x_2),\xi(x_3)]_\h+\xi[x_1,x_2,x_3]_\g.
\end{eqnarray}
\end{pro}
\pf Let $(\mathfrak{g}\oplus \mathfrak{h},[\cdot,\cdot,\cdot]_{(\rho^{1},\nu^1,\omega^{1})})$ and $(\mathfrak{g}\oplus \mathfrak{h},[\cdot,\cdot,\cdot]_{(\rho^{2},\nu^2,\omega^{2})})$ be  two non-abelian extensions of $\g$
by $\frkh$. Assume that the two extensions are isomorphic. Then there is a $3$-Lie algebra morphism $\theta:\mathfrak{g}\oplus \mathfrak{h}\lon \mathfrak{g}\oplus \mathfrak{h}$, such that we have the following commutative diagram:
\[\begin{CD}
0@>>>\mathfrak{h}@>\iota >>\mathfrak{g}\oplus \mathfrak{h}_{(\rho^{2},\nu^2,\omega^{2})}@>\pr>>{\mathfrak{g}}  @>>>0\\
@.    @|                       @V\theta VV                                                   @|               @.\\
0@>>>\mathfrak{h}@>\iota >>\mathfrak{g}\oplus \mathfrak{h}_{(\rho^{1},\nu^1,\omega^{1})}@>\pr>>{\mathfrak{g}}  @>>>0,
\end{CD}\]
where $\iota$ is the inclusion and $\pr$ is the projection. Since for all $x\in\g$, $\pr(\theta(x))=x$, we can assume that     $\theta(x+u)=x-\xi(x)+u$ for some linear map $\xi:\mathfrak{g}\lon \mathfrak{h}$. By
$$\theta[x_1,v_2,v_3]_{(\rho^{2},\nu^2,\omega^{2})}=[\theta(x_1),\theta(v_2),\theta(v_3)]_{(\rho^{1},\nu^1,\omega^{1})},$$
we can deduce that \eqref{iso1} holds. By
$$\theta[x_1,x_2,v_3]_{(\rho^{2},\nu^2,\omega^{2})}=[\theta(x_1),\theta(x_2),\theta(v_3)]_{(\rho^{1},\nu^1,\omega^{1})},$$
we can deduce that \eqref{iso2} holds. By
$$\theta[x_1,x_2,x_3]_{(\rho^{2},\nu^2,\omega^{2})}=[\theta(x_1),\theta(x_2),\theta(x_3)]_{(\rho^{1},\nu^1,\omega^{1})},$$
we can deduce that \eqref{iso3} holds.\qed

\section{Non-abelian extensions in terms of Maurer Cartan elements}
In \cite{NR bracket of n-Lie}, the author constructed a graded Lie algebra structure by which one can describe an $n$-Leibniz algebra structure as a canonical structure. Here, we give the precise formulas for the 3-Lie algebra case.

  Set $C^p(\g,\g)=\Hom(\wedge^2\frkg\otimes\stackrel{(p)}{\cdots }\otimes\wedge^2\frkg\wedge\frkg,\g)$ and $C^*(\g,\g)=\oplus_{p}C^{p}(\g,\g)$. Let $\alpha\in C^p(\g,\g),\beta\in C^q(\g,\g),\, p,q\geq 0$. Let $\frkX_i=x_i\wedge y_i\in \wedge^2\g$ for $i=1,2,\cdots ,p+q$ and $x_i,y_i\in\g$. A permutation $\sigma\in\mathbb S_n$ is called an $(i,n-i)$-unshuffle if $\sigma(1)<\cdots <\sigma(i)$ and $\sigma(i+1)<\cdots <\sigma(n)$. If $i=0$ or $n$, we assume $\sigma=\Id$. The set of all $(i,n-i)$-unshuffles will be denoted by $unsh(i,n-i)$.

  \begin{thm}{\rm (\cite{NR bracket of n-Lie})}\label{thm:gradelie}
  The graded vector space $C^*(\g,\g)$ equipped with the graded commutator bracket
\begin{eqnarray}
\nrn{\alpha,\beta}= \alpha\circ \beta-(-1)^{pq}\beta\circ \alpha,
\end{eqnarray}
is a graded Lie algebra where $\alpha\circ \beta\in C^{p+q}(\g,\g)$ is defined by
\begin{eqnarray*}
&&\alpha\circ \beta(\frkX_1,\cdots ,\frkX_{p+q},x)=
\sum_{k=0}^{p-1}(-1)^{kq}\Big(\sum_{\sigma\in unsh(k,q)}\\&&(-1)^{\sigma}\big(\alpha(\frkX_{\sigma(1)},\cdots ,\frkX_{\sigma(k)},\beta(\frkX_{\sigma(k+1)},\cdots ,\frkX_{\sigma(k+q)},x_{k+q+1})\wedge y_{k+q+1},\frkX_{k+q+2},\cdots ,\frkX_{p+q},x)\\
&&+\alpha(\frkX_{\sigma(1)},\cdots ,\frkX_{\sigma(k)},x_{k+q+1}\wedge\beta(\frkX_{\sigma(k+1)},\cdots ,\frkX_{\sigma(k+q)},y_{k+q+1}) ,\frkX_{k+q+2},\cdots ,\frkX_{p+q},x)\big)\Big)\\
&&+\sum_{\sigma\in unsh(p,q)}(-1)^{pq}(-1)^{\sigma}\alpha(\frkX_{\sigma(1)},\cdots ,\frkX_{\sigma(p)},\beta(\frkX_{\sigma(p)},\cdots ,\frkX_{\sigma(p+q-1)}, \frkX_{\sigma(p+q)},x)).
\end{eqnarray*}
\end{thm}

Furthermore, $(C^*(\g,\g),\nrn{\cdot,\cdot},\overline{\delta})$ is a DGLA, where $\overline{\delta}$ is given by $\overline{\delta}P=(-1)^p\delta P$ for all $P\in C^{p}(\g,\g)$, and $\delta $ is the coboundary operator of $\g$ with  coefficients in the adjoint representation. See \cite{Liu-Makhlouf-Sheng}  for more details.
\begin{rmk}
The coboundary operator $\delta$ associated to the adjoint representation of the $3$-Lie algebra $\g$ can be written as $\delta P=(-1)^p\nrn{\mu_\g,P}$, for all $P\in C^{p}(\g,\g)$,   where $\mu_{\g}\in C^{1}(\g,\g)$ is the $3$-Lie algebra structure on $\g$, i.e. $\mu_\g(x,y,z)=[x,y,z]_{\g}$. Thus, we have $\overline{\delta}P=\nrn{\mu_\g,P}.$
\end{rmk}

Now, we describe non-abelian extensions using Maurer-Cartan elements. The set $MC(L)$ of {\bf Maurer-Cartan elements} of a DGLA $(L,[\cdot,\cdot],\dM)$ is defined by
$$MC(L)\triangleq \{P\in L_1\mid\dM P+\frac{1}{2}[P,P]=0\}.$$
Moreover, $P_0,P_1\in MC(L)$ are called {\bf gauge equivalent} if and only if there exists an element $\xi\in L_0$ such that
\begin{eqnarray}
P_1=e^{\ad_\xi}P_0-\frac{e^{\ad_\xi}-1}{\ad_\xi}\dM\xi.
\end{eqnarray}
We can define a path between $P_0$ and $P_1$. Let
\begin{eqnarray*}
P(t)=e^{t\ad_\xi}P_0-\frac{e^{t\ad_\xi}-1}{\ad_\xi}\dM\xi.
\end{eqnarray*}
Then $P(t)$ is a power series of $t$ in $MC(L)$. We have $P(0)=P_0$ and $P(1)=P_1$. The set of the gauge equivalence classes of $MC(L)$ is denoted by $\huaMC(L)$.

Let $(\g,[\cdot,\cdot,\cdot]_\g)$ and $(\frkh,[\cdot,\cdot,\cdot]_{\frkh})$ be two 3-Lie algebras. Let $\g\oplus \h$ be the 3-Lie algebra direct sum of $\g$ and $\h$, where the bracket is defined by
$$[x+u,y+v,z+w]=\mu_\g(x,y,z)+\mu_\h(u,v,w).$$
Then there is a DGLA   $(C(\g\oplus \h,\g\oplus \h), \nrn{\cdot,\cdot},\overline{\delta})$. Define $C^k_>(\g\oplus\h,\h)\subset C^k(\g\oplus\h,\h)$ by
$$C^k(\g\oplus\h,\h)=C^k_>(\g\oplus\h,\h)\oplus C^k(\h,\h).$$
Denote by $C_>(\g\oplus\h,\h)=\oplus_{k}C^k_>(\g\oplus\h,\h)$, which is a graded vector space.

\begin{lem}\label{lem:dgla}
  With the above notations, we have $\overline{\delta} (C^k_>(\g\oplus\h,\h))\subset C^{k+1}_>(\g\oplus\h,\h)$, and $({C}_>(\g\oplus\h,\h),\nrn{\cdot,\cdot},\overline{\delta})$ is a sub-DGLA of $(C(\g\oplus\h,\g\oplus\h), \nrn{\cdot,\cdot},\overline{\delta})$.    Furthermore, its degree $0$ part $C^0_>(\g\oplus\h,\h)=\Hom(\g,\h)$ is abelian.
\end{lem}

\pf By the definition of the bracket $\nrn{\cdot,\cdot}$ and $\overline{\delta}$, we obtain that $(C(\g\oplus \h, \h),\nrn{\cdot,\cdot},\overline{\delta})$ is a sub-DGLA of $(C(\g\oplus \h,\g\oplus \h), \nrn{\cdot,\cdot},\overline{\delta})$. For $\alpha\in C^p_>(\g\oplus\h,\h)$, we can regard it as $\alpha\in C^p(\g\oplus\h,\h)$ such that $\alpha|_{C^p(\h,\h)}=0$. Moreover, for $\alpha\in C^p_>(\g\oplus\h,\h),\beta\in C^q_>(\g\oplus\h,\h)$, we have $\nrn{\alpha,\beta}|_{C^{p+q}(\h,\h)}=0$ and $\overline{\delta}(\alpha)|_{C^{p+1}(\h,\h)}=\nrn{\mu_\g+\mu_\h,\alpha}|_{C^{p+1}(\h,\h)}=0$. Thus, we obtain $(C(\g\oplus \h, \h),\nrn{\cdot,\cdot},\overline{\delta})$ is a sub-DGLA of $(C(\g\oplus \h,\h), \nrn{\cdot,\cdot},\overline{\delta})$. Therefore, $({C}_>(\g\oplus\h,\h),\nrn{\cdot,\cdot},\overline{\delta})$ is a sub-DGLA of $(C(\g\oplus\h,\g\oplus\h), \nrn{\cdot,\cdot},\overline{\delta})$. Obviously, $C^0_>(\g\oplus\h,\h)=\Hom(\g,\h)$ is abelian.\qed



\begin{pro}\label{pro:JMC}
The following two statements are equivalent:
\begin{itemize}
  \item[\rm(a)] $(\g\oplus \h,[\cdot,\cdot,\cdot]_{(\rho,\nu,\omega)})$ is a $3$-Lie algebra, which is a non-abelian extension of $\g$ by $\h$;

  \item[\rm(b)] $\rho+\nu+\omega$ is a Maurer-Cartan element of the DGLA $(C_>(\g\oplus \h, \h),\nrn{\cdot,\cdot},\overline{\delta})$.
\end{itemize}
\end{pro}

\pf By  Proposition \ref{3-lie}, $(\g\oplus\h,[\cdot,\cdot,\cdot]_{(\rho,\nu,\omega)})$ is a $3$-Lie algebra if and only if Eqs. \eqref{p1}-\eqref{p11} hold.

If $c=\rho+\nu+\omega$ is a Maurer-Cartan element, we have
$$(\overline{\delta}c+\frac{1}{2}\nrn{c,c})(e_1\wedge e_2,e_3\wedge e_4,e_5)=0, \quad \forall e_i=x_i+v_i\in\g\oplus \h.$$
By straightforward computations, we have
\begin{eqnarray*}
\overline{\delta}c(e_1\wedge e_2,e_3\wedge e_4,e_5)&=&-(\delta c)(e_1\wedge e_2,e_3\wedge e_4,e_5)\\
                                                   &=&c([e_1,e_2,e_3]\wedge e_4,e_5)+c(e_3\wedge [e_1,e_2,e_4],e_5)\\
                                                   &&+c(e_3\wedge e_4,[e_1,e_2,e_5])-c(e_1\wedge e_2,[e_3,e_4,e_5])\\
                                                   &&-[e_1,e_2,c(e_3\wedge e_4,e_5)]+[e_3,e_4,c(e_1\wedge e_2,e_5)]\\
                                                   &&+[e_4,e_5,c(e_1\wedge e_2,e_3)]+[e_5,e_3,c(e_1\wedge e_2,e_4)].
\end{eqnarray*}
Recall the definition of  the bracket $[\cdot,\cdot,\cdot]$ and the $c$, we have
\begin{eqnarray*}
[e_1,e_2,e_3] &=&[x_1,x_2,x_3]_{\frkg}+[v_1,v_2,v_3]_{\frkh},\\
c(e_1,e_2,e_3)&=&\rho(x_1,x_2)v_3+\rho(x_2,x_3)v_1+\rho(x_3,x_1)v_2\\&&+\nu(x_1)(v_2,v_3)+\nu(x_2)(v_3,v_1)+\nu(x_3)(v_1,v_2)
+\omega(x_1,x_2,x_3).
\end{eqnarray*}

Moreover, we have
\begin{eqnarray*}
c([e_1,e_2,e_3]\wedge e_4,e_5)&=&\rho([x_1,x_2,x_3]_\g,x_4)(v_5)+\rho(x_4,x_5)([v_1,v_2,v_3]_{\frkh})+\rho(x_5,[x_1,x_2,x_3]_\g)(v_4)\\
                               &&+\nu([x_1,x_2,x_3]_\g)(v_4,v_5)+\nu(x_4)(v_5,[v_1,v_2,v_3]_{\frkh})+\nu(x_5)([v_1,v_2,v_3]_{\frkh},v_4)\\
                               &&+\omega([x_1,x_2,x_3]_{\frkg},x_4,x_5),\\
c(e_3\wedge [e_1,e_2,e_4],e_5)&=&\rho(x_3,[x_1,x_2,x_4]_\g)(v_5)+\rho([x_1,x_2,x_4]_\g,x_5)(v_3)+\rho(x_5,x_3)([v_1,v_2,v_4]_{\frkh})\\
                               &&+\nu(x_3)([v_1,v_2,v_4]_{\frkh},v_5)+\nu([x_1,x_2,x_4]_\g)(v_5,v_3)+\nu(x_5)(v_3,[v_1,v_2,v_4]_{\frkh})\\
                               &&+\omega(x_3,[x_1,x_2,x_4]_\g,x_5),\\
c(e_3\wedge e_4,[e_1,e_2,e_5])&=&\rho(x_3,x_4)([v_1,v_2,v_5]_{\frkh})+\rho(x_4,[x_1,x_2,x_5]_\g)(v_3)+\rho([x_1,x_2,x_5]_\g,x_3)(v_4)\\
                               &&+\nu(x_3)(v_4,[v_1,v_2,v_5]_{\frkh})+\nu(x_4)([v_1,v_2,v_5]_{\frkh},v_3)+\nu([x_1,x_2,x_5]_\g)(v_3,v_4)\\
                               &&+\omega(x_3,x_4,[x_1,x_2,x_5]_\g),\\
c(e_1\wedge e_2,[e_3,e_4,e_5])&=&\rho(x_1,x_2)([v_3,v_4,v_5]_{\frkh})+\rho(x_2,[x_3,x_4,x_5]_\g)(v_1)+\rho([x_3,x_4,x_5]_\g,x_1)(v_2)\\
                               &&+\nu(x_1)(v_2,[v_3,v_4,v_5]_{\frkh})+\nu(x_2)([v_3,v_4,v_5]_{\frkh},v_1)+\nu([x_3,x_4,x_5]_\g)(v_1,v_2)\\
                               &&+\omega(x_1,x_2,[x_3,x_4,x_5]_\g),\\
\label{}[e_1,e_2,c(e_3\wedge e_4,e_5)]&=&[v_1,v_2,\rho(x_3,x_4)(v_5)]_\h+[v_1,v_2,\rho(x_4,x_5)(v_3)]_\h+[v_1,v_2,\rho(x_5,x_3)(v_4)]_\h\\
                                      &&+[v_1,v_2,\nu(x_3)(v_4,v_5)]_\h+[v_1,v_2,\nu(x_4)(v_5,v_3)]_\h+[v_1,v_2,\nu(x_5)(v_3,v_4)]_\h\\
                                      &&+[v_1,v_2,\omega(x_3,x_4,x_5)]_\h,\\
\label{}[e_3,e_4,c(e_1\wedge e_2,e_5)]&=&[v_3,v_4,\rho(x_1,x_2)(v_5)]_\h+[v_3,v_4,\rho(x_2,x_5)(v_1)]_\h+[v_3,v_4,\rho(x_5,x_1)(v_2)]_\h\\
                                      &&+[v_3,v_4,\nu(x_1)(v_2,v_5)]_\h+[v_3,v_4,\nu(x_2)(v_5,v_1)]_\h+[v_3,v_4,\nu(x_5)(v_1,v_2)]_\h\\
                                      &&+[v_3,v_4,\omega(x_1,x_2,x_5)]_\h,\\
\label{}[e_4,e_5,c(e_1\wedge e_2,e_3)]&=&[v_4,v_5,\rho(x_1,x_2)(v_3)]_\h+[v_4,v_5,\rho(x_2,x_3)(v_1)]_\h+[v_4,v_5,\rho(x_3,x_1)(v_2)]_\h\\
                                      &&+[v_4,v_5,\nu(x_1)(v_2,v_3)]_\h+[v_4,v_5,\nu(x_2)(v_3,v_1)]_\h+[v_4,v_5,\nu(x_3)(v_1,v_2)]_\h\\
                                      &&+[v_4,v_5,\omega(x_1,x_2,x_3)]_\h,\\
\label{}[e_5,e_3,c(e_1\wedge e_2,e_4)]&=&[v_5,v_3,\rho(x_1,x_2)(v_4)]_\h+[v_5,v_3,\rho(x_2,x_4)(v_1)]_\h+[v_5,v_3,\rho(x_4,x_1)(v_2)]_\h\\
                                      &&+[v_5,v_3,\nu(x_1)(v_2,v_4)]_\h+[v_5,v_3,\nu(x_2)(v_4,v_1)]_\h+[v_5,v_3,\nu(x_4)(v_1,v_2)]_\h\\
                                      &&+[v_5,v_3,\omega(x_1,x_2,x_4)]_\h.
\end{eqnarray*}

\emptycomment{
\begin{eqnarray*}
c([e_1,e_2,e_3]\wedge e_4,e_5)&=&\rho([x_1,x_2,x_3]_\g,x_4)(v_5)+\rho(x_4,x_5)(\omega(x_1,x_2,x_3))+\rho(x_4,x_5)(\rho(x_1,x_2)v_3)\\
                               &&+\rho(x_4,x_5)(\rho(x_2,x_3)v_1)+\rho(x_4,x_5)(\rho(x_3,x_1)v_2)+\rho(x_4,x_5)(\nu(x_1)(v_2,v_3))\\
                               &&+\rho(x_4,x_5)(\nu(x_2)(v_3,v_1))+\rho(x_4,x_5)(\nu(x_3)(v_1,v_2))+\rho(x_4,x_5)([v_1,v_2,v_3]_{\frkh})\\
                               &&+\rho(x_5,[x_1,x_2,x_3]_\g)(v_4)+\nu([x_1,x_2,x_3]_\g)(v_4,v_5)+\nu(x_4)(v_5,\omega(x_1,x_2,x_3))\\
                               &&+\nu(x_4)(v_5,\rho(x_1,x_2)v_3)+\nu(x_4)(v_5,\rho(x_2,x_3)v_1)+\nu(x_4)(v_5,\rho(x_3,x_1)v_2)\\
                               &&+\nu(x_4)(v_5,\nu(x_1)(v_2,v_3))+\nu(x_4)(v_5,\nu(x_2)(v_3,v_1))+\nu(x_4)(v_5,\nu(x_3)(v_1,v_2))\\
                               &&+\nu(x_4)(v_5,[v_1,v_2,v_3]_{\frkh})+\nu(x_5)(\omega(x_1,x_2,x_3),v_4)+\nu(x_5)(\rho(x_1,x_2)v_3,v_4)\\
                               &&+\nu(x_5)(\rho(x_2,x_3)v_1,v_4)+\nu(x_5)(\rho(x_3,x_1)v_2,v_4)+\nu(x_5)(\nu(x_1)(v_2,v_3),v_4)\\
                               &&+\nu(x_5)(\nu(x_2)(v_3,v_1),v_4)+\nu(x_5)(\nu(x_3)(v_1,v_2),v_4)+\nu(x_5)([v_1,v_2,v_3]_{\frkh},v_4)\\
                               &&+\omega([x_1,x_2,x_3]_{\frkg},x_4,x_5).
\end{eqnarray*}
}

Furthermore, by the definition of the bracket in Theorem \ref{thm:gradelie}, we have
\begin{eqnarray*}
(\frac{1}{2}\nrn{c,c})(e_1\wedge e_2,e_3\wedge e_4,e_5)&=&(c\circ c)(e_1\wedge e_2,e_3\wedge e_4,e_5)\\
                                                       &=&c(c(e_1\wedge e_2,e_3)\wedge e_4,e_5)+c(e_3\wedge c(e_1\wedge e_2,e_4),e_5)\\
                                                        &&-c(e_1\wedge e_2,c(e_3\wedge e_4,e_5))+c(e_3\wedge e_4,c(e_1\wedge e_2,e_5)).
\end{eqnarray*}
Moreover, we have
\begin{eqnarray*}
c(c(e_1\wedge e_2,e_3)\wedge e_4,e_5)&=&\rho(x_4,x_5)(\omega(x_1,x_2,x_3))+\rho(x_4,x_5)(\rho(x_1,x_2)v_3)+\rho(x_4,x_5)(\rho(x_2,x_3)v_1)\\
                               &&+\rho(x_4,x_5)(\rho(x_3,x_1)v_2)+\rho(x_4,x_5)(\nu(x_1)(v_2,v_3))+\rho(x_4,x_5)(\nu(x_2)(v_3,v_1))\\
                               &&+\rho(x_4,x_5)(\nu(x_3)(v_1,v_2))+\nu(x_4)(v_5,\omega(x_1,x_2,x_3))+\nu(x_4)(v_5,\rho(x_1,x_2)v_3)\\
                               &&+\nu(x_4)(v_5,\rho(x_2,x_3)v_1)+\nu(x_4)(v_5,\rho(x_3,x_1)v_2)+\nu(x_4)(v_5,\nu(x_1)(v_2,v_3))\\
                               &&+\nu(x_4)(v_5,\nu(x_2)(v_3,v_1))+\nu(x_4)(v_5,\nu(x_3)(v_1,v_2))+\nu(x_5)(\omega(x_1,x_2,x_3),v_4)\\
                               &&+\nu(x_5)(\rho(x_1,x_2)v_3,v_4)+\nu(x_5)(\rho(x_2,x_3)v_1,v_4)+\nu(x_5)(\rho(x_3,x_1)v_2,v_4)\\
                               &&+\nu(x_5)(\nu(x_1)(v_2,v_3),v_4)+\nu(x_5)(\nu(x_2)(v_3,v_1),v_4)+\nu(x_5)(\nu(x_3)(v_1,v_2),v_4),\\
c(e_3\wedge c(e_1\wedge e_2,e_4),e_5)&=&\rho(x_5,x_3)(\rho(x_1,x_2)v_4)+\rho(x_5,x_3)(\rho(x_2,x_4)v_1)+\rho(x_5,x_3)(\rho(x_4,x_1)v_2)\\
                                      &&+\rho(x_5,x_3)(\nu(x_1)(v_2,v_4))+\rho(x_5,x_3)(\nu(x_2)(v_4,v_1))+\rho(x_5,x_3)(\nu(x_4)(v_1,v_2))\\
                                      &&+\rho(x_5,x_3)(\omega(x_1,x_2,x_4))+\nu(x_3)(\rho(x_1,x_2)v_4,v_5)+\nu(x_3)(\rho(x_2,x_4)v_1,v_5)\\
                                      &&+\nu(x_3)(\rho(x_4,x_1)v_2,v_5)+\nu(x_3)(\nu(x_1)(v_2,v_4),v_5)+\nu(x_3)(\nu(x_2)(v_4,v_1),v_5)\\
                                      &&+\nu(x_3)(\nu(x_4)(v_1,v_2),v_5)+\nu(x_3)(\omega(x_1,x_2,x_4),v_5)+\nu(x_5)(v_3,\rho(x_1,x_2)v_4)\\
                                      &&+\nu(x_5)(v_3,\rho(x_2,x_4)v_1)+\nu(x_5)(v_3,\rho(x_4,x_1)v_2)+\nu(x_5)(v_3,\nu(x_1)(v_2,v_4))\\
                                      &&+\nu(x_5)(v_3,\nu(x_2)(v_4,v_1))+\nu(x_5)(v_3,\nu(x_4)(v_1,v_2))+\nu(x_5)(v_3,\omega(x_1,x_2,x_4)),\\
c(e_1\wedge e_2,c(e_3\wedge e_4,e_5))&=&\rho(x_1,x_2)(\rho(x_3,x_4)v_5)+\rho(x_1,x_2)(\rho(x_4,x_5)v_3)+\rho(x_1,x_2)(\rho(x_5,x_3)v_4)\\
                                      &&+\rho(x_1,x_2)(\nu(x_3)(v_4,v_5))+\rho(x_1,x_2)(\nu(x_4)(v_5,v_3))+\rho(x_1,x_2)(\nu(x_5)(v_3,v_4))\\
                                      &&+\rho(x_1,x_2)(\omega(x_3,x_4,x_5))+\nu(x_1)(v_2,\rho(x_3,x_4)v_5)+\nu(x_1)(v_2,\rho(x_4,x_5)v_3)\\
                                      &&+\nu(x_1)(v_2,\rho(x_5,x_3)v_4)+\nu(x_1)(v_2,\nu(x_3)(v_4,v_5))+\nu(x_1)(v_2,\nu(x_4)(v_5,v_3))\\
                                      &&+\nu(x_1)(v_2,\nu(x_5)(v_3,v_4))+\nu(x_1)(v_2,\omega(x_3,x_4,x_5))+\nu(x_2)(\rho(x_3,x_4)v_5,v_1)\\
                                      &&+\nu(x_2)(\rho(x_4,x_5)v_3,v_1)+\nu(x_2)(\rho(x_5,x_3)v_4,v_1)+\nu(x_2)(\nu(x_3)(v_4,v_5),v_1)\\
                                      &&+\nu(x_2)(\nu(x_4)(v_5,v_3),v_1)+\nu(x_2)(\nu(x_5)(v_3,v_4),v_1)+\nu(x_2)(\omega(x_3,x_4,x_5),v_1),\\
c(e_3\wedge e_4,c(e_1\wedge e_2,e_5))&=&\rho(x_3,x_4)(\rho(x_1,x_2)v_5)+\rho(x_3,x_4)(\rho(x_2,x_5)v_1)+\rho(x_3,x_4)(\rho(x_5,x_1)v_2)\\
                                      &&+\rho(x_3,x_4)(\nu(x_1)(v_2,v_5))+\rho(x_3,x_4)(\nu(x_2)(v_5,v_1))+\rho(x_3,x_4)(\nu(x_5)(v_1,v_2))\\
                                      &&+\rho(x_3,x_4)(\omega(x_1,x_2,x_5))+\nu(x_3)(v_4,\rho(x_1,x_2)v_5)+\nu(x_3)(v_4,\rho(x_2,x_5)v_1)\\
                                      &&+\nu(x_3)(v_4,\rho(x_5,x_1)v_2)+\nu(x_3)(v_4,\nu(x_1)(v_2,v_5))+\nu(x_3)(v_4,\nu(x_2)(v_5,v_1))\\
                                      &&+\nu(x_3)(v_4,\nu(x_5)(v_1,v_2))+\nu(x_3)(v_4,\omega(x_1,x_2,x_5))+\nu(x_4)(\rho(x_1,x_2)v_5,v_3)\\
                                      &&+\nu(x_4)(\rho(x_2,x_5)v_1,v_3)+\nu(x_4)(\rho(x_5,x_1)v_2,v_3)+\nu(x_4)(\nu(x_1)(v_2,v_5),v_3)\\
                                      &&+\nu(x_4)(\nu(x_2)(v_5,v_1),v_3)+\nu(x_4)(\nu(x_5)(v_1,v_2),v_3)+\nu(x_4)(\omega(x_1,x_2,x_5),v_3).
\end{eqnarray*}
Thus, $c=\rho+\nu+\omega$ is a Maurer-Cartan element if and only if \eqref{p1}-\eqref{p11} hold.\qed

\begin{cor}\label{MCtoext1}
Let $\g$ and $\h$ be two $3$-Lie algebras. Then there is a one-to-one correspondence between  non-abelian extensions of the $3$-Lie algebra $\g$ by $\h$ and Maurer-Cartan elements in the DGLA $(C_>(\g\oplus \h, \h),\nrn{\cdot,\cdot},\overline{\delta})$.
\end{cor}

\begin{thm}\label{MCtoext2}
Let $\g$ and $\h$ be two $3$-Lie algebras. Then the isomorphism classes of non-abelian extensions $\g$ by $\h$ one-to-one correspond to the gauge equivalence classes of Maurer-Cartan elements in the DGLA $(C_>(\g\oplus \h,\h),\nrn{\cdot,\cdot},\overline{\delta})$.
\end{thm}

\pf Two elements $c=\rho+\nu+\omega$ and $c'=\rho'+\nu'+\omega'$ in $MC(L)$ are equivalent if there exists $\xi\in \Hom(\g,\h)$ such that
$$c'=e^{\ad_\xi}c-\frac{e^{\ad_\xi}-1}{\ad_\xi}\overline{\delta}\xi.$$
More precisely, for all $e_i=x_i+v_i\in\g\oplus \h$, we have
\begin{eqnarray*}
c'(e_1\wedge e_2,e_3)&=&\big(({\Id}+\ad_{\xi}+\frac{1}{2!}\ad_{\xi}^2+\frac{1}{3!}\ad_{\xi}^3+\cdots+\frac{1}{n!}\ad_{\xi}^{n}+\cdots)c\big)(e_1\wedge e_2,e_3)\\
                      &&-\big(({\Id}+\frac{1}{2!}\ad_{\xi}+\frac{1}{3!}\ad_{\xi}^2+\cdots+\frac{1}{n!}\ad_{\xi}^{n-1}+\cdots)\overline{\delta} \xi\big)(e_1\wedge e_2,e_3).
\end{eqnarray*}
Furthermore, by the bracket in Theorem \ref{thm:gradelie}, we have
\begin{eqnarray*}
\nrn{\xi,c}(e_1\wedge e_2,e_3)&=&(\xi\circ c)(e_1\wedge e_2,e_3)-(c\circ\xi)(e_1\wedge e_2,e_3)\\
                              &=&-\big(c(\xi(e_1)\wedge e_2,e_3)+c(e_1\wedge \xi(e_2),e_3)+c(e_1\wedge e_2,\xi(e_3))\big)\\
                              &=&-c(\xi(x_1)\wedge e_2,e_3)-c(e_1\wedge \xi(x_2),e_3)-c(e_1\wedge e_2,\xi(x_3))\\
                              &=&-\rho(x_2,x_3)\xi(x_1)-\nu(x_2)(v_3,\xi(x_1))-\nu(x_3)(\xi(x_1),v_2)\\
                              &&-\rho(x_3,x_1)\xi(x_2)-\nu(x_1)(\xi(x_2),v_3)-\nu(x_3)(v_1,\xi(x_2))\\
                              &&-\rho(x_1,x_2)\xi(x_3)-\nu(x_1)(v_2,\xi(x_3))-\nu(x_2)(\xi(x_3),v_1).
\end{eqnarray*}
Thus, we have
\begin{eqnarray*}
&&\nrn{\xi,\nrn{\xi,c}}(e_1\wedge e_2,e_3)\\
&=&-\nrn{\xi,c}(\xi(x_1)\wedge e_2,e_3)-\nrn{\xi,c}(e_1\wedge \xi(x_2),e_3)-\nrn{\xi,c}(e_1\wedge e_2,\xi(x_3))\\
&=&2\nu(x_1)(\xi(x_2),\xi(x_3))+2\nu(x_2)(\xi(x_3),\xi(x_1))+2\nu(x_3)(\xi(x_1),\xi(x_2)).
\end{eqnarray*}
Moreover, we have
\begin{eqnarray*}
\nrn{\xi,\nrn{\xi,\nrn{\xi,c}}}(e_1\wedge e_2,e_3)=0.
\end{eqnarray*}
More generally, for $n\ge 3$
$$\ad_{\xi}^n c=0.$$
For all $e_i=x_i+v_i\in\g\oplus \h$, we have
\begin{eqnarray*}
\overline{\delta}\xi(e_1\wedge e_2,e_3)&=&\nrn{\mu_\g+\mu_\h,\xi}(e_1\wedge e_2,e_3)\\
                         &=&((\mu_\g+\mu_\h)\circ\xi)(e_1\wedge e_2,e_3)-(\xi\circ(\mu_\g+\mu_\h))(e_1\wedge e_2,e_3)\\
                         &=&(\mu_\g+\mu_\h)(\xi(e_1)\wedge e_2,e_3)+(\mu_\g+\mu_\h)(e_1\wedge \xi(e_2),e_3)\\
                         &&+(\mu_\g+\mu_\h)(e_1\wedge e_2,\xi(e_3))-\xi((\mu_\g+\mu_\h)(e_1\wedge e_2,e_3))\\
                         &=&[\xi(x_1),v_2,v_3]_\h+[v_1,\xi(x_2),v_3]_\h+[v_1,v_2,\xi(x_3)]_\h-\xi[x_1,x_2,x_3]_\g.
\end{eqnarray*}
Thus, we have
\begin{eqnarray*}
\nrn{\xi,\overline{\delta}\xi}(e_1\wedge e_2,e_3)&=&-\overline{\delta}\xi(\xi(x_1)\wedge e_2,e_3)-\overline{\delta}\xi(e_1\wedge \xi(x_2),e_3)-\overline{\delta}\xi(e_1\wedge e_2,\xi(x_3))\\
                                   &=&-2[\xi(x_1),\xi(x_2),v_3]_\h-2[v_1,\xi(x_2),\xi(x_3)]_\h-2[\xi(x_1),v_2,\xi(x_3)]_\h,
\end{eqnarray*}
and
\begin{eqnarray*}
&&\nrn{\xi,\nrn{\xi,\overline{\delta}\xi}}(e_1\wedge e_2,e_3)\\&=&-\nrn{\xi,\overline{\delta}\xi}(\xi(x_1)\wedge e_2,e_3)-\nrn{\xi,\overline{\delta}\xi}(e_1\wedge \xi(x_2),e_3)-\nrn{\xi,\overline{\delta}\xi}(e_1\wedge e_2,\xi(x_3))\\
                                             &=&6[\xi(x_1),\xi(x_2),\xi(x_3)]_\h.
\end{eqnarray*}
Moreover, we have
$$
\nrn{\xi,\nrn{\xi,\nrn{\xi,\overline{\delta}\xi}}}(e_1\wedge e_2,e_3)=0.
$$
More generally, for $n\ge 3$
$$\ad_{\xi}^n \overline{\delta}\xi=0.$$
Therefore,we have
\begin{eqnarray*}
c'=(c+\nrn{\xi,c}+\frac{1}{2!}\nrn{\xi,\nrn{\xi,c}})-(\overline{\delta}\xi+\frac{1}{2!}\nrn{\xi,\overline{\delta}\xi}+\frac{1}{3!}\nrn{\xi,\nrn{\xi,\overline{\delta}\xi}}).
\end{eqnarray*}
Thus, two elements $c=\rho+\nu+\omega$ and $c'=\rho'+\nu'+\omega'$ in $MC(L)$ are equivalent if and only \eqref{iso1}-\eqref{iso3} hold. \qed

\section{Non-abelian extensions of Leibniz algebras}
In this section, we always assume that $(\mathfrak{g}\oplus \mathfrak{h},[\cdot,\cdot]_{(\rho,\nu,\omega)})$ is a non-abelian extension of the $3$-Lie algebra $\g$ by $\h$. We aim to  analyze the corresponding Leibniz algebra structure on the space of fundamental
objects. Note that $\wedge^2(\g\oplus \h)\cong((\wedge^2\h)\oplus(\g\otimes\h))\oplus(\wedge^2\g)$ naturally. We use $[\cdot,\cdot]_{\rm \hat{F}}$ to denote the Leibniz bracket on the space of fundamental objects of the 3-Lie algebra $(\mathfrak{g}\oplus \mathfrak{h},[\cdot,\cdot]_{(\rho,\nu,\omega)})$.

First we introduce a Leibniz algebra structure on $(\wedge^2\h)\oplus(\g\otimes\h)$. Define a linear map $\{\cdot,\cdot\}:((\wedge^2\h)\oplus(\g\otimes\h))\otimes((\wedge^2\h)\oplus(\g\otimes\h))\lon(\wedge^2\h)\oplus(\g\otimes\h)$ by
\begin{eqnarray}
\nonumber&&\{u_1\wedge v_1+x_1\otimes w_1,u_2\wedge v_2+x_2\otimes w_2\}\\\nonumber
&&=[u_1,v_1,u_2]_\h\wedge v_2+u_2\wedge [u_1,v_1,v_2]_\h+\nu(x_2)(u_1,v_1)\wedge w_2+x_2\otimes [u_1,v_1,w_2]_\h\\
                                                           &&+\nu(x_1)(w_1,u_2)\wedge v_2+u_2\wedge \nu(x_1)(w_1,v_2)
                                                           -\rho(x_1,x_2)(w_1)\wedge w_2+x_2\otimes\nu(x_1)(w_1,w_2).
\end{eqnarray}

\begin{pro}
With the above notations, $((\wedge^2\h)\oplus(\g\otimes\h),\{\cdot,\cdot\})$ is a Leibniz algebra.
\end{pro}

\pf By direct computation, we have
\begin{eqnarray*}
\{u_1\wedge v_1+x_1\otimes w_1,u_2\wedge v_2+x_2\otimes w_2\}=[u_1\wedge v_1+x_1\otimes w_1,u_2\wedge v_2+x_2\otimes w_2]_{\rm \hat{F}}.
\end{eqnarray*}
Thus, $((\wedge^2\h)\oplus(\g\otimes\h),\{\cdot,\cdot\})$ is a Leibniz subalgebra of the Leibniz algebra $(\wedge^2(\g\oplus \h),[\cdot,\cdot]_{\rm \hat{F}})$.  \qed\vspace{3mm}

We define $\varpi:(\wedge^2\g)\otimes(\wedge^2\g)\lon(\wedge^2\h)\oplus(\g\otimes\h),~l:(\wedge^2\g)\lon\gl((\wedge^2\h)\oplus(\g\otimes\h))$, and $r:(\wedge^2\g)\lon\gl((\wedge^2\h)\oplus(\g\otimes\h))$ respectively by
\begin{eqnarray}
\varpi(x\wedge y,z\wedge t)&=&- t\otimes \omega(x,y,z)+z\otimes\omega(x,y,t),\\
\nonumber l(x\wedge y)(u\wedge v+z\wedge w)&=&\rho(x,y)(u)\wedge v+u\wedge\rho(x,y)(v)+[x,y,z]_\g\otimes w\\
                                    &&+\omega(x,y,z)\wedge w+z\otimes\rho(x,y)(w),\\
\nonumber r(x\wedge y)(u\wedge v+z\wedge w)&=&-y\otimes \nu(x)(u,v)+x\otimes\nu(y)(u,v)\\
                                    &&  y\otimes \rho(z,x)(w)-x\otimes\rho(z,y)(w),
\end{eqnarray}
for all $x,y,z,t\in\g,~u,v,w\in\h.$

Now we are ready to give the main result of this section.
\begin{thm}
Let $(\g,[\cdot,\cdot,\cdot]_\g)$ and $(\h,[\cdot,\cdot,\cdot]_\h)$ be two $3$-Lie algebras and $(\mathfrak{g}\oplus \mathfrak{h},[\cdot,\cdot]_{(\rho,\nu,\omega)})$  a non-abelian extension of the $3$-Lie algebra $\g$ by $\h$. Then the Leibniz algebra $(\wedge^2(\g\oplus \h),[\cdot,\cdot]_{\rm \hat{F}})$ is a non-abelian extension of the Leibniz algebra $(\wedge^2\g,[\cdot,\cdot]_{\rm F})$ by the Leibniz algebra $((\wedge^2\h)\oplus(\g\otimes\h),\{\cdot,\cdot\})$.
\end{thm}

\pf One can show that conditions \eqref{l Der}-\eqref{cocycle} in Proposition \ref{non-abelian extension of LB} hold directly. Thus, $(\wedge^2(\g\oplus \h),[\cdot,\cdot]_{\rm \hat{F}})$ is a non-abelian extension of the Leibniz algebra $(\wedge^2\g,[\cdot,\cdot]_{\rm F})$ by the Leibniz algebra $((\wedge^2\h)\oplus(\g\otimes\h),\{\cdot,\cdot\})$. Here we use a different approach to prove this theorem. Using the isomorphism between $\wedge^2(\g\oplus \h)$ and $((\wedge^2\h)\oplus(\g\otimes\h))\oplus(\wedge^2\g)$, the Leibniz algebra structure on $\wedge^2(\g\oplus \h)$ is given by
\begin{eqnarray*}
&&[u_1\wedge v_1+x_1\otimes w_1+y_1\wedge z_1,u_2\wedge v_2+x_2\otimes w_2+y_2\wedge z_2]_{\rm \hat{F}}\\
&&=[u_1\wedge v_1+x_1\otimes w_1,u_2\wedge v_2+x_2\otimes w_2]_{\rm \hat{F}}+[y_1\wedge z_1,u_2\wedge v_2+x_2\otimes w_2]_{\rm \hat{F}}\\
&&+[u_1\wedge v_1+x_1\otimes w_1,y_2\wedge z_2]_{\rm \hat{F}}+[y_1\wedge z_1,y_2\wedge z_2]_{\rm \hat{F}}\\
&&=\{u_1\wedge v_1+x_1\otimes w_1,u_2\wedge v_2+x_2\otimes w_2\}+l(y_1\wedge z_1)(u_2\wedge v_2+x_2\otimes w_2)\\
&&+r(y_2\wedge z_2)(u_1\wedge v_1+x_1\otimes w_1)+\varpi(y_1\wedge z_1,y_2\wedge z_2)+[y_1\wedge z_1,y_2\wedge z_2]_{\rm F}.
\end{eqnarray*}
Thus, by \eqref{bracket of nbe}, we deduce that $(\wedge^2(\g\oplus \h),[\cdot,\cdot]_{\rm \hat{F}})$ is a non-abelian extension of the Leibniz algebra $(\wedge^2\g,[\cdot,\cdot]_{\rm F})$ by the Leibniz algebra $((\wedge^2\h)\oplus(\g\otimes\h),\{\cdot,\cdot\})$.  \qed

\end{document}